\documentclass[12pt,journal,onecolumn]{IEEEtran}
\usepackage{amsmath, epsfig, cite}
\usepackage{amssymb}
\usepackage{amsfonts}
\usepackage{latexsym}
\usepackage{longtable}
\usepackage{tabularx}
\usepackage{color}
\usepackage{multirow}
\usepackage{pgf}
\usepackage{tikz}
\usepackage{float}
\usetikzlibrary{arrows, automata, positioning, calc}
\usepackage{array}
\usepackage{colortbl}
\usepackage{comment}
\usepackage{dblfloatfix}
\usepackage{hyperref}
\usepackage{graphicx}
\usepackage{makecell}
\usepackage{diagbox}

\newcommand{\E}{\mathbb{E}}
\newcommand{\F}{\mathbb{F}}

\newcommand{\cC}{\mathcal C}
\newcommand{\cX}{\mathcal X}
\newcommand{\tr}{\textup{Tr}}
\newcommand{\Aut}{\textup{Aut}}

\newtheorem{thm}{Theorem}

\newtheorem{cor}[thm]{Corollary}

\newtheorem{lemma}[thm]{Lemma}

\newtheorem{example}{Example}
\newtheorem{remark}{Remark}

\definecolor{burntorange}{rgb}{0.8, 0.33, 0.0}

\begin{document}

\title{\textbf{A construction of minimal \\
linear codes from partial difference sets}}


\author{
{ Ran~Tao, \hspace{0.1cm}
Tao~Feng, \hspace{0.1cm}
Weicong~Li$^*$}

\thanks{Ran Tao and Tao Feng are with the School of Mathematical Sciences, Zhejiang University, Hangzhou 310027, Zhejiang P.R. China, (e-mail: \{rant,tfeng\}@zju.edu.cn). Weicong Li is with the Department of Mathematics, Southern University of Science and Technology,  ShenZhen 518055, Guangdong P.R. China, (e-mail: liwc3@sustech.edu.cn).}
\thanks{$^*$ Corresponding author: Weicong Li}}
\maketitle

\begin{abstract}
In this paper, we study a class of linear codes defined by characteristic functions of certain subsets of a finite field. We derive a sufficient and necessary condition for such a code to be a minimal linear code by a character-theoretical approach. We obtain new three-weight or four-weight minimal linear codes that do not satisfy the Ashikhmin-Barg condition by using partial difference sets. We show that our construction yields minimal linear codes that do not arise from cutting vectorial blocking sets, and also discuss their applications in secret sharing schemes.   \vspace*{10pt}
\end{abstract}


\section{Introduction}
\label{sec:intro}

Let $q$ be a prime power and $\mathbb{F}_{q}$ be the finite field with $q$ elements. Let $\mathcal{C}$ be an $[n,k]$ linear code over $\mathbb{F}_{q}$, that is, a $k$-dimensional subspace of $\mathbb{F}_{q}^{n}$. The \textit{support} of a codeword $c=(c_{1},c_{2},\ldots,c_{n})\in\mathcal{C}$ is defined by $\text{supp}(c)=\{1\le i\le n:\,c_{i}\neq0\}$. The \textit{Hamming weight} of a codeword $c$ is $wt(c):=|\text{supp}(c)|$, the size of $\text{supp}(c)$. For any two codewords $u,\,v\in\mathcal{C}$, we say that $v$ \textit{covers} $u$ and write $u\preceq v$ if $\text{supp}(u)\subseteq\text{supp}(v)$. Clearly, $au\preceq v$ for all $a\in\mathbb{F}_{q}$ if $u\preceq v$. A codeword $c\in\mathcal{C}$ is called \textit{minimal} if $c$ covers only the codewords $\lambda c$ with $\lambda\in\mathbb{F}_{q}$, and no other codewords in $\mathcal{C}$. The linear code $\mathcal{C}$ is called \textit{minimal} if every codeword $c\in\mathcal{C}$ is minimal.

As a special type of linear codes, minimal linear codes have many applications in sharing schemes \cite{carlet2005linear, massey1993minimal, massey1995some, yuan2005secret} and secure two-party computation \cite{chabanne2013towards}.

One useful and easy criterion was given by Ashikhmin and Barg \cite{AB1998} for a linear code $\cC$ to be minimal, which is referred to as AB condition in the literature \cite{bonini2019minimal}. \vspace*{10pt}
\begin{lemma}\cite{AB1998}
\label{ABbound}
Let $\cC$ be a linear code over $\F_q$. Set $\omega_{\text{min}}$ and $\omega_{\text{max}}$ to be the minimum and maximum nonzero weights of the code $\cC$ respectively. Then $\cC$ is minimal if
\begin{equation}\label{eq_ABcond}
\frac{\omega_{\text{min}}}{\omega_{\text{max}}}>\frac{q-1}{q}.
\end{equation}
\end{lemma}
This condition is sufficient but not necessary for linear codes to be minimal. A sufficient and necessary condition was given by Heng et al. in \cite{heng2018minimal}:
\begin{lemma}\cite{heng2018minimal}
\label{SNCHeng}
 A linear code $\cC\subseteq \F_{q}^{n}$ over $\F_q$ is minimal if and only if
 \[\sum_{a\in\F_{q}^*}wt(c'+ac)\ne (q-1)wt(c')-wt(c).\]
 for any two $\F_{q}$-linearly independent codewords $c,c'$ of $\cC$.
\end{lemma}

In \cite{heng2018minimal}, the authors applied Lemma \ref{SNCHeng} to construct minimal ternary linear codes that do not satisfy the AB condition. Up to now, many minimal linear codes not satisfying the AB condition have been constructed. We summarize some known minimal linear codes in Table \ref{table_codes} below. These minimal linear codes do not satisfy the AB condition when their parameters satisfy certain conditions. For more constructions and information on minimal linear codes, please refer to \cite{bartoli2019minimal, bonini2019minimal, ding2018minimal, heng2018minimal, mesnager2019minimal, Shi2020, xu2019three, Xu2020, zhang2019four}.

\begin{table}[!htbp]
  \centering
  \caption{Some known minimal linear codes}
  \label{table_codes}
  \scalebox{0.75}{
    \begin{tabular}{|c|c|c|c|c|c|}
      \hline
      $[n,\,k]$ & minimal distance $d$ & The numbers of weights &  Notes& Approach &References\\ \hline
      $[3^m-1,\,m+1]$ & $\sum_{j=1}^k 2^j \binom{m}{j} $ &  $\le m+2$ & $m\geq 5,\, 2\leq k\leq \lfloor (m-1)/2 \rfloor$ &  Boolean functions& \cite{heng2018minimal} \\ \hline
      \multirow{5}*{$[2^m-1,\, m+1]$} & \multirow{2}*{$\min({s(2^{t}-1),\,2^{m-1}-s})$} & \multirow{2}*{$3$ or $4$}  &  $m\geq  6 $ even, $t=m/2$ &\multirow{5}*{Boolean functions} &\multirow{5}*{ \cite{ding2018minimal}}\\

      & & & $s \notin\{1,2^t,2^t+1\}$ & & \\
      \cline{2-4}

       &\multirow{2}*{$2^{m-1}-2^{m-s-1}(s-1)$}  &$s+3$($s$ odd) & \multirow{2}*{$m\ge 7$ odd, $s=(m+1)/2$} &  &\\
     &  & $s+2$($s$ even) & & &  \\ \cline{2-4}
    & $\sum_{j=1}^k\binom{m}{j}$ &$\le m+2$ & $m\ge7$,\, $2\leq k\leq \lfloor (m-3)/2\rfloor$ & &\\
      \hline
      $[p^m-1,\,m]$& \multirow{2}*{$(p-1)^2p^{m-2}$} & $3$ & \multirow{3}*{$m >2 $} & \multirow{3}*{$p$-ary functions}& \multirow{4}*{\cite{xu2019three}} \\ \cline{1-1}\cline{3-3}
    $[p^m-1,\,m-1]$ & &\multirow{2}*{2}& & &\\ \cline{1-2}
       $[p^m-1,\,m]$ &  $p^{m-1}(p-2)$ & & & & \\  \cline{1-5}
         $[p^m-1,\,m+1]$ & - & $3$ or $4$ &$m=2t, t\ge 2$ & partial spreads &  \\ \hline
      \multirow{3}*{$[p^m-1,\,m+1]$} & $p^{m-s-1}(p-1)(s(p-1)+1)$  & $\leq s+4$ & {$(p-1)(p^{s-2}-s)>1$} & \multirow{3}*{Maiorana-McFarland functions} &\multirow{3}*{\cite{Xu2020}} \\ \cline{2-4}
      & \multirow{2}*{$a(p^{m-s}-p^{m-s-1})(p^k-1)$} & \multirow{2}*{$6$} & $k\geq 2$ and $(k,p)\ne (2,2)$  && \\
      && &$s=2k$,$2\leq a\leq (p-1)p^{k-2}$ & &\\ \hline
      \multirow{2}*{ $[p^m-1,\, m+1]$} &\multirow{2}*{-} & {4 or 5} &{Theorem 3.8}& Characteristic functions& \multirow{2}*{\cite{mesnager2019minimal}} \\ \cline{3-4}
      && 6 or 7  & Theorem 3.12 &corresponding to subspaces &  \\ \hline

      $[q^m-1,\, m+1]$ & -  & - &  - &Cutting blocking sets & \cite{bartoli2019minimal}, \cite{bonini2019minimal} \\ \hline

  \end{tabular}}
\end{table}

In \cite{bonini2019minimal}, Bonini and Borello developed a method to construct minimal linear codes by using cutting blocking set. An affine $k$-\textit{blocking set} is a subset of an $n$-dimensional affine space intersecting all $(n-k)$-dimensional affine subspaces. An affine 1-blocking set is also called an \textit{affine blocking set}. A \textit{vectorial $k$-blocking set} is a subset of an $n$-dimensional affine space not containing the origin and intersecting all $(n-k)$-dimensional affine subspaces through the origin. A vectorial 1-blocking set is also called a \textit{vectorial blocking set}. A \textit{vectorial $(k,s)$-blocking set} is a vectorial $k$-blocking set that does not contain an $s$-dimensional affine subspace through the origin. A vectorial $k$-blocking set is \textit{cutting} if its intersection with every $(n-k)$-dimensional affine  subspace through the origin is not contained in any other $(n-k)$-dimensional affine subspace through the origin.\vspace*{10pt}

\begin{lemma}
\cite[Theorem 4.6]{bonini2019minimal}
\label{lem_cutting}
Let $f:\,\F_{q}^m\rightarrow\F_{q}$ be a function that is not linear. If
\begin{enumerate}
	\item[1)]$V(f)^{*}=\{{\bf x}\in\F_{q}^{m}\setminus\{{0}\}:\, f({\bf x})=0\}$ is an $m$-dimensional cutting vectorial $(1,m-1)$-blocking set in the affine space $\F_{q}^{m}$; 
	\item[2)] for every nonzero vector $\bf{v}$, if exists $\bf{x}$ such that $f(\bf{x})+ v\cdot x=0$ and $f(\bf{x})$ is different from  $0$, 
\end{enumerate} then  $\cC(f)=\{\left(uf(\bf{x})+\bf{v}\cdot \bf{x}\right)_{{\bf x}\in\F_{q}^{m}\setminus\{0\}}: \,u\in\F_{q},\,{\bf v}\in\F_{q}^{m}\}$ is a $[q^{m}-1,m+1]$ minimal linear codes over $\F_{q}$.
\end{lemma}

Let $m$ be a positive integer and $\F_{q^m}^{*}$ be the multiplicative group of finite field $\F_{q^m}$. Take a proper subset $D$ of $\F_{q^m}^{*}$ and define its characteristic function as
\[f_{D}(x)=\left\{\begin{array}{ll}
1, &\text{if } x \in D,\\
0, & \text{if } x \in \F_{q^m}^{*}\setminus D.
\end{array}\right.\]
We define a linear code
\begin{equation}
\label{eqn_CfD}
\cC(f_D)=\left\{\left(uf_{D}(x)+\tr_{\F_{q^m}/\F_q}(vx)\right)_{x\in \F_{q^m}^*}: (u,v) \in V\right\},
\end{equation}
where $\tr_{\F_{q^m}/\F_q}$ is the trace function from $\F_{q^m}$ to $\F_{q}$.

Here is a brief summary of the main results of this paper. We give a necessary and sufficient condition on $D$ such that the code $\cC(f_{D})$ is minimal in Theorem \ref{thm_SNC2}. By taking $D$ to be an $\F_{q}^{*}$-invariant partial difference set, we deduce sufficient conditions on the parameters of $D$ such that $\cC(f_{D})$ is minimal in Theorem \ref{thm_mainpds}. We also show that if the parameters of $D$ satisfy certain conditions then the code $\cC(f_D)$ does not satisfy the AB condition. Our construction yields minimal linear codes that do not arise from cutting vectorial blocking sets, cf. Examples \ref{example_cyc}, \ref{example_sporadic}. In the case $D$ is an $\F_{q}^{*}$-invariant partial difference set, we determine the weight distribution of $\cC(f_{D})$. In Section \ref{subsec_auto}, we show that each automorphism of the subset $D$ induces an automorphism of the code $\cC(f_{D})$. In particular, we obtain a minimal linear code $\cC(f_{D})$ with a large automorphism group in some cases, which potentially has a fast decoding algorithm.

This paper is organized as follows. In Section \ref{sec:Preliminaries}, we present some preliminary results on polynomials, characters, strongly regular graphs and partial difference sets. In Section \ref{sec:Mainresults}, we give our construction of minimal linear codes from partial difference sets and study their properties. In Section \ref{sec:Applications}, we present the applications of our minimal linear codes in secret sharing schemes. Finally, we conclude this paper in Section \ref{sec:Concluding}.

\section{Preliminaries}
\label{sec:Preliminaries}
Throughout this paper, let $q$ be a prime power and $m$ be a positive integer. Let $\F_{q}$ and $\F_{q^m}$ be the finite fields with $q$ elements and $q^m$ elements respectively. Let $\F_{q}^{*}$ (resp. $\F_{q^m}^{*}$) be the multiplicative group of nonzero elements of $\F_{q}$ (resp. $\F_{q^m}$).  For any subset $S$ of $\F_{q^m}$, we write  $\langle S\rangle$ for the $\F_q$-linear subspace spanned by $S$, and we say that $S$ is \textit{$\F_{q}^{*}$-invariant} if $\lambda s\in S$ for each $\lambda\in\F_{q}^{*}$ and $s\in S$. For an extension field $\E$ of the finite field $\F$ with finite degree $[\E:\F]$, we use $\tr_{\E/\F}$ to denote the trace function from $\E$ to $\F$.

\subsection{Some basic facts on finite fields}\
\label{sec:basicfacts}
Let $G$ be a finite abelian group of order $n$. A \textit{character} $\phi$ of $(G,+)$ is a group homomorphism $\phi:\,G\rightarrow \mathbb{C}^{\times}$, where $\mathbb{C}^{\times}$ is the multiplicative group of nonzero complex numbers. The \textit{ trivial} character $\phi_{0}$ of $G$ is defined by $\phi_{0}(g)=1$ for all $g\in G$. The set $\widehat{G}$ of characters of $G$ forms an abelian group with identity $\phi_{0}$, which is isomorphic to $G$, cf. \cite[Chapter 5]{LidlFF}. The group $\widehat{G}$ is called the character group of $G$.

An \textit{additive character} $\psi$ of $\mathbb{F}_{q}$ is a character of the additive group $\left(\mathbb{F}_{q},+\right)$. For each $a\in\mathbb{F}_{q}$, we can define  an additive character $\psi_{a}$ of $\F_q$ via
\[\psi_{a}(x)=\zeta_{p}^{\text{Tr}_{\F_{q}/\F_{p}}(ax)},\]
where $\zeta_{p}=e^{\frac{2\pi i}{p}}$ is a primitive $p$-th root of unity. The character $\psi_{1}$ is called the \textit{canonical additive character}, which we denote by $\psi$. The character group of $(\F_{q},+)$ is given by $\widehat{(\F_{q},+)}=\{\psi_a:\,a \in \F_{q}\}$. For an extension $\F_{q^m}$ of $\F_q$, its canonical additive character $\Psi$ can be obtained from  $\psi$, i.e., \[\Psi(x)=\psi(\tr_{\F_{q^m}/\F_{q}}(x)).\]
The character group of $(\F_{q^m},+)$ is given by $\widehat{(\F_{q^m},+)}=\{\Psi_a:\,a \in \F_{q^m}\}$, where $\Psi_a(x)=\Psi(ax)$ for all $x\in\F_{q^m}$. There is an important property about the additive characters, cf. \cite{LidlFF}:
\begin{align}
\label{eqn_charsum}
\sum_{\lambda\in\F_{q}}\Psi(\lambda x)&=\sum_{\lambda\in\F_{q}}\psi(\lambda\tr_{\F_{q^m}/\F_{q}}(x))
=\left\{\begin{array}{ll}
                                      q, & \text{ if }\tr_{\F_{q^m}/\F_{q}}(x)=0,\\
                                      0, & \text{ otherwise.}
                                    \end{array}
\right.
\end{align}
For each $a\in\F_{q^{m}}$ and a subset $A\subseteq\F_{q^m}$, we define
\begin{equation}
\label{eqn_psiaA}
 \Psi_a(A):=\sum_{x\in A} \Psi(ax).
\end{equation}

We regard $\F_{q^m}$ as an $\F_q$-linear vector space of dimension $m$. For each $a\in\F_{q^m}^{*}$, we define
\[L(a):=\{x\in\F_{q^m}:\tr_{\F_{q^m}/\F_{q}}(xa)=0\},\]
which is an $\F_q$-linear subspace of $\F_{q^m}$ of dimension $m-1$. For a subset $S$ of $\F_{q^m}^*$, we define
\begin{equation}\label{eqn_LS}
L(S):=\{x\in\F_{q^m}: \tr_{\F_{q^m}/\F_{q}}(xs)=0 \text{ for all  }s\in S\}.
\end{equation}
Since $L(S)$ is the intersection of the $L(s)$'s with $s\in S$, it is also an $\F_q$-linear subspace of $\F_{q^m}$.\vspace*{10pt}
\begin{lemma}\label{lem_1}
Let $L$ be as defined in Eqn. \eqref{eqn_LS}. Let $V_1$ be an $\F_{q}$-linear subspace of $\F_{q^m}$ and $U_1$ be a subset of $V_1$. Then $\langle U_1\rangle=V_1$ if and only if $L(U_1) \subseteq L(V_1)$, where $
\langle U_{1}\rangle$ is the $\F_q$-linear subspace spanned by $U_{1}$.
\end{lemma}
\begin{IEEEproof}
	For any subset $S\subseteq\F_{q^m}$, it is routine to check that $L(\langle S\rangle)=L(S)$ and $L(L(S))=\langle S\rangle$ by comparing their dimensions. Therefore, $V_{1}\subseteq\langle U_{1}\rangle$ if and only if $L(U_1)=L(\langle  U_{1}\rangle)\subseteq L(V_1)$, and the claim follows.
\end{IEEEproof}\vspace*{10pt}
\begin{lemma}\label{lem_2}
	Let $a\in\F_{q^m}$ and $A \subseteq \F_{q^m}$. Then $A \subseteq L(a)$ if and only if $\Psi_{\lambda a}(A)=|A|$ for all $\lambda\in\F_q$, where $\Psi_{\lambda a}(A)$ is defined in Eqn. \eqref{eqn_psiaA}. In particular, if $A$ is $\F_q^*$-invariant, then $L(A)=\{a\in\F_{q^m}:\,\Psi_a(A)=|A|\}.$
\end{lemma}
\begin{IEEEproof}
Suppose that $A\subseteq L(a)$. Then $\tr_{\F_{q^m}/\F_{q}}(xa)=0 $ for all $x\in A$, and thus $\Psi_{\lambda a}(x)=\psi\left(\lambda\tr_{\F_{q^m}/\F_{q}}( ax)\right)$ $=1$ for all $x\in A$ and $\lambda\in\F_q$. It follows that $\Psi_{\lambda a}(A)=|A|$ for all $\lambda\in\F_q$ by definition of $\Psi_{\lambda a}(A)$ as in Eqn. \eqref{eqn_psiaA}. Conversely, suppose that $\Psi_{\lambda a}(A)=|A|$ for all $\lambda\in\F_q$, we compute that
	\begin{align*}
	q|A|&=\sum_{\lambda \in\F_q} \Psi_{\lambda a}(A)=\sum_{x\in A}\sum_{\lambda \in\F_q}\Psi(\lambda ax)\\
	&=\sum_{x\in A}\sum_{\lambda \in\F_q}\psi\left(\lambda \tr_{\F_{q^m}/\F_{q}}(ax)\right).
	\end{align*}
	By Eqn. \eqref{eqn_charsum}, we deduce that $\tr_{\F_{q^m}/\F_q}(ax)=0$ for all $x \in A$, i.e., $A\subseteq L(a)$. In the case $A$ is $\F_q^*$-invariant, we have $\Psi_{\lambda a}(A)=\Psi_a(\lambda A)=\Psi_a(A)$ for any $\lambda\in \F_q^*$. In addition, it is obvious that $\Psi_{0}(A)=|A|$. It follows that
\begin{align*}
 L(A)&=\{a\in\F_{q^m}: \tr_{\F_{q^m}/\F_{q}}(xa)=0 \text{ for all  }x\in A\}\\
 &= \{a\in \F_{q^m}:\,  A\subseteq L(a)\}\\
 &=\{a\in \F_{q^m}: \Psi_{a}(A)=|A|\}.
 \end{align*}
	This completes the proof.
\end{IEEEproof}\vspace*{10pt}

A polynomial $f(X)\in\F_{q^m}[X]$ is called a permutation polynomial if the associated function $f:\,a\mapsto f(a)$ from $\F_{q^m}$ to itself is a permutation. A polynomial of the form $f(X)=\sum_{i=0}^{n}a_iX^{q^i}$ with $a_i\in\F_{q^m}$ is called a $q$-polynomial over $\F_{q^m}$. It is called \textit{reduced} if $n\le m-1$. There is a one to one correspondence between the reduced $q$-polynomials over $\F_{q^m}$ and $\F_q$-linear transformations of $\F_{q^m}$, cf. \cite{LidlFF}. Let $f(X)=\sum_{i=0}^{m-1}a_iX^{q^i}$ be a reduced $q$-polynomial over $\F_{q^m}$. The \textit{trace dual} of $f(X)$ is the (unique) reduced $q$-polynomial $\tilde{f}(X)$ such that $\tr_{\F_{q^m}/\F_{q}}(f(x)y)=\tr_{\F_{q^m}/\F_q}(\tilde{f}(y)x)$ for all $x,\,y\in\F_{q^m}$. A direct computation shows that
\begin{equation}\label{tracedual}
\tilde{f}(X)=\sum\limits_{i=0}^{m-1}a_{m-i}^{q^i}X^{q^i}.
\end{equation}

\subsection{Strongly regular graphs and partial difference sets}\
\label{sec:pds}
A strongly regular graph (srg) $(v,k,\lambda,\mu)$ is a simple and undirected graph $\Gamma$, neither complete nor edgeless, that has the following properties:
\begin{enumerate}
\item It is a $k$-regular graph of order $v$.

\item Every two adjacent vertices have $\lambda$ common neighbors.
\item Every two nonadjacent vertices have $\mu$ common neighbors.

\end{enumerate}
An eigenvalue is called restricted if it has an eigenvector perpendicular to the all-one vector, cf. \cite{Brouwer2012Spec}. The following result is well known. \vspace*{10pt}

\begin{lemma}\cite[Section 10.2]{Godsil2001Alg}\label{lem_eigenvalue}
	Let $\Gamma$ be an srg $(v,k,\lambda,\mu)$ with restricted eigenvalues $\theta_1$ and $\theta_2$, where $\theta_1>0>\theta_2$. Then \[\theta_{1}=\frac{(\lambda-\mu)+\sqrt{(\lambda-\mu)^2+4(k-\mu)}}{2},\] \[\theta_{2}=\frac{(\lambda-\mu)-\sqrt{(\lambda-\mu)^2+4(k-\mu)}}{2}\]
	with multiplicities \[m_{1}=\frac{1}{2}\left((v-1)-\frac{2k+(v-1)(\lambda-\mu)}{\sqrt{(\lambda-\mu)^2+4(k-\mu)}}\right),\] \[m_{2}=\frac{1}{2}\left((v-1)+\frac{2k+(v-1)(\lambda-\mu)}{\sqrt{(\lambda-\mu)^2+4(k-\mu)}}\right)\] respectively.
\end{lemma}

Let $G$ be an (additive) abelian group of order $v$, and let $D$ be a subset of $G$ such that $0\notin D$ and $-D=D$, where $-D=\{-d:\,\,d\in D\}$. The \textit{Cayley graph} $\Gamma=\text{Cay}(G,D)$ with the connection set $D$ is the graph whose vertices are the elements of $G$ such that $x\sim y$ if and only if $x-y\in D$ for any $x,y\in G$. In the case that $\text{Cay}(G,D)$ is an srg $(v,k,\lambda,\mu)$, the connection set $D$ is called a \textit{partial difference set} (PDS for short) with parameters $(v,k,\lambda,\mu)$. For $\Gamma=\text{Cay}(G,D)$ with $G$ abelian,  $\{\phi(D):=\sum_{d\in D}\phi(d),\, \phi\in\widehat{G}\}$ consists of all the eigenvalues of $\Gamma$, where $\widehat{G}$ is the character group of $G$. For a PDS $D$ with parameters $(v,k,\lambda,\mu)$, we have $k=\phi_{0}(D)=|D|$, where $\phi_{0}$ is the trivial character of $G$. Please refer to the survey \cite{Ma1994} for more details on PDS.

An srg $(v,k,\lambda,\mu)$ is of \textit{Latin square} type (resp. \textit{negative Latin square} type) if there exist positive integers $n,\,r$ such that
\[(v,k,\lambda,\mu)=(n^2,r(n-\epsilon),\epsilon n+r^2-3\epsilon r,r^2-\epsilon r)\]
where $\epsilon=1$ (resp. $\epsilon=-1$). Correspondingly, a PDS $D$ is called a Latin square type (resp. negative Latin square type) PDS if the srg $\text{Cay}(G,D)$ is of the corresponding type.

In \cite{brouwer1999cyclotomy}, the authors gave a construction of partial difference sets from cyclotomy, which we describe now. Set $q=p^e$ with $p$ prime. Let $m$ be a positive integer such that 2 divides $em$ and $\gamma$ be a fixed primitive element of $\mathbb{F}_{q^m}$. For a proper divisor $N$ of $q^m-1$, we define the $N$-th \textit{cyclotomic classes} of $\mathbb{F}_{q^m}$ as $C_{i}=\{\gamma^{jN+i}:\,0\leq j\leq \frac{q^m-1}{N}-1\}$, where $0\leq i \leq N-1$. The set $C_{0}$ is a subgroup of $\mathbb{F}_{q^m}^{*}$ of index $N$, and $C_{i}=\gamma^i C_{0}$ for $0\leq i\leq N-1$.\medskip

\begin{lemma}\cite{brouwer1999cyclotomy}
\label{PDS}
Take notation as above and suppose that $N$ is a proper divisor of $q^m-1$ such that $N\neq1$ and $p^{\ell_{1}}\equiv-1(\text{mod }N)$ for some positive integer $\ell_{1}$. Choose $\ell_{1}$ minimal and write $em=2\ell_{1}t$. Take a proper subset $J\subset\mathbb{Z}_{N}$ of size $u$. If $q$ is odd, we further assume that $N|\frac{q^m-1}{2}$ and $J+\frac{q^m-1}{2}\equiv J(\text{mod }N)$. Set
\begin{equation*}
\label{DefDJ}
D=D_{J}=\bigcup_{j\in J}C_{j}.
\end{equation*}
Then the graph $\text{Cay}(\mathbb{F}_{q^m},D)$ is strongly regular with eigenvalues
\[k=|D|=\frac{q^m-1}{N}u,\  \text{with\ multiplicity\ 1};\]
\[\theta_{1}=\frac{u}{N}(-1+(-1)^t \sqrt{q^{m}}),\  \text{with\ multiplicity\ } q^m-1-k;\]
\[\theta_{2}=\theta_{1}+(-1)^{t+1} \sqrt{q^{m}},\  \text{with\ multiplicity\ } k.\]
To be specific, for $i=0,\,1,\,\ldots,\,N-1$, we have
\[\Psi(\gamma^i D)=\left\{\begin{array}{lll}\theta_2, & \text{if\  } \varepsilon^t=1, i\in-J(\text{mod\ }N)\text{\ or\ }\\
&\varepsilon^t=-1, i\in-J+N/2(\text{mod\ }N), \\ \theta_{1}, &  \text{\ otherwise,} \end{array}\right.\]
where $\varepsilon=\left\{\begin{array}{ll}-1, & \text{ if }N\text{\ is\ even\ and\ }\frac{p^{\ell_{1}}+1}{N}\text{\ is\ odd},  \\ 1, &  \text{\ otherwise}. \end{array}\right.$

The graph $\text{Cay}(\mathbb{F}_{q^m},D)$ is of Latin square type (resp. negative Latin square type) if $t$ is odd (resp. even).
\end{lemma}\vspace*{10pt}

The following lemma is an easy consequence of the fact that $\mathbb{F}_{q}^{*}$ is the subgroup of $\mathbb{F}_{q^m}^{*}$ consisting of all nonzero $\frac{q^m-1}{q-1}$-th powers.\vspace*{10pt}

\begin{lemma}
\label{lem_Fqinv}
Take notation as above, and let $D=D_{J}=\bigcup_{j\in J}C_{j}$ be as defined in Lemma \ref{PDS}. Then $D$ is $\mathbb{F}_{q}^{*}$-invariant if and only if the set $J$ is invariant under the map $\rho:\,j\rightarrow j+\frac{q^m-1}{q-1}(\text{mod }N)$.
\end{lemma}

\section{Minimal linear codes from partial difference sets}
\label{sec:Mainresults}
Set $V=\F_q\times \F_{q^m}$, and view it as a vector space of dimension $m+1$ over $\F_q$. Let  $B$ be a bilinear form on $V$ such that
\begin{equation}
\label{Bilinear}
B((u,v),\,(x,y))= ux+\tr_{\F_{q^m}/\F_q}(vy)
\end{equation}
for $(u,v),\,(x,y)\in V$. For $(a,b)\in V$, the \textit{perp} of $(a,b)$ is defined by $(a,b)^{\perp}=\{(x,y)\in V: B((a,b),\,(x,y))=0 \}$. Correspondingly for any subset $S$ of $V$, we have
\[S^{\perp}=\{(x,y)\in V:\, B((a,b),\,(x,y))=0,\, \forall\, (a,b)\in S \}. \]
Take a subset $M$ of $V$, we can construct a linear code from $M$ as follows:

\begin{equation}\label{eqn_codeM}
\cC(M):=\{c(u,v)=(B((u,v),\,(a_{i},b_{i})))_{1\leq i\leq n}: (u,v) \in V \},
\end{equation}
where $M=\{(a_{1},b_{1}),\,(a_{2},b_{2}),\ldots,(a_{n},b_{n})\}$ is called the \textit{defining set} of $\cC(M)$.

Let $\cC(f_D)$ be the linear code as defined in Eqn. \eqref{eqn_CfD}. It is clear that
 \[\cC(f_D)=\{(B((u,v),\,(f_{D}(x),x)))_{x\in \F_{q^m}^*}: (u,v) \in V \}.\]
By the definition of $B$ as in Eqn. (\ref{Bilinear}), the code $\cC(f_D)$ can be constructed from the defining set. We define the following subset of $V$:
\begin{equation}\label{eqn_MD}
M_{D}:=\{(1,r):r\in D\} \cup\{(0,r): r \in \overline{D}\},
\end{equation}
where $\overline{D}=\F_{q^m}^* \setminus D$. It is easy to check that $\cC(f_D)=\cC(M_D)$.
The dimension of $\cC(f_D)$ $\left(\text{i.e., }\cC(M_{D})\right)$ has been determined if $f_D$ is not linear.\vspace*{10pt}
\begin{lemma}\cite{bonini2019minimal}\label{lem_code_para}
Let $D$ be a proper subset of $\F_{q^m}^*$ such that the characteristic function $f_D$ is not linear. Then the linear code $\cC(f_D)$ defined in Eqn. \eqref{eqn_CfD} has length $q^m-1$ and dimension $m+1$ over $\F_q$.
\end{lemma}

\begin{lemma}
\label{lem_fnonlinear}
Let $D$ be a proper $\F_{q}^{*}$-invariant subset of $\F_{q^m}^*$ and $f_{D}$ be the characteristic function defined in Eqn. \eqref{eqn_CfD}. If $q>2$, then $f_{D}$ is not linear.
\end{lemma}

\begin{IEEEproof}
Suppose that $f_{D}$ is linear, i.e., $f_{D}(x+y)=f_{D}(x)+f_{D}(y)$ for any $x,y\in\F_{q^{m}}^{*}$. We take any $x\in D$ and choose $\lambda\in\F_q$ such that $\lambda\not\in\{0,-1\}$. Since $D$ is $\F_q^*$-invariant, both $\lambda x$ and $(\lambda+1)x$ are in $D$. We deduce from $f_{D}((\lambda+1)x)=f_{D}(\lambda x)+f_{D}(x)$ that $1+1=1$ in $\F_q$, which is impossible. This completes the proof.
\end{IEEEproof}

\subsection{A sufficient and necessary condition for $\cC(M_{D})$ to be minimal}\

For a subset $D$ of $\F_{q^m}^*$, set $\overline{D}=\F_{q^m}^* \setminus D$. Let $L(S)$ be defined in Eqn. \eqref{eqn_LS} for a subset $S$ of $\F_{q^m}^*$. For $y\in \F_q$, $z\in \F_{q^m}^*$, define
\begin{equation}
\label{eqn_overlineDz}
\overline{D}_z:=\{ x\in \overline{D}:\tr_{\F_{q^m}/\F_{q}}(xz)=0 \},\\
\end{equation}
\begin{equation}
\label{eqn_Dyz}
D_{(y,z)}:=\{ x\in D:\tr_{\F_{q^m}/\F_q}(xz)=-y\},
\end{equation}
\begin{equation}\label{eqn_pyz}
P_{(y,z)}:=L\left(D_{(y,z)}D_{(y,z)}^{(-1)}\cup \overline{D}_z\right),
\end{equation}
 where $D_{(y,z)}D_{(y,z)}^{(-1)}=\{d_i-d_j:\,d_i,d_j\in D_{(y,z)}\}$.

Recall that for a subset $S$ of $\F_{q^m}$, $\langle S\rangle$ is an $\F_q$-linear subspace spanned by $S$. We can apply Theorem 3.2 and Theorem 3.3 of \cite{lu2019parameters} to get the following theorem directly. Here we briefly repeat the proof in our language.

\begin{thm}\label{thm_SNC1}
Take notation as above. Suppose that $D$ is a subset of $\F_{q^m}^*$ and $M_D$ is defined in  Eqn. (\ref{eqn_MD}). Let $\cC(M_{D})$ be a linear code defined as in Eqn. \eqref{eqn_codeM} with the defining set $M_D$. For each $(y,z)\in V$ with $(y,z) \ne (0,0)$, the codeword $c(y,z)\in\cC(M_{D})$ is minimal if and only if $\langle(y,z)^{\perp}\cap M_D\rangle=(y,z)^{\perp}$. In particular, $\cC(M_{D})$ is minimal if and only if $\langle(y,z)^{\perp}\cap M_D\rangle=(y,z)^{\perp}$ for all $(y,z)\in V\setminus \{(0,0)\}$. 
\end{thm}
\begin{IEEEproof}
For any two codewords $c(y_{1},z_{1}),\,c(y_{2},z_{2})\in \cC(M_{D})$, we have that $c(y_{1},z_{1})\preceq c(y_{2},z_{2})$, i.e.,
$\text{supp}(c(y_{1},z_{1}))\subseteq\text{supp}(c(y_{2},z_{2}))$, if and only if $(y_{2},z_{2})^{\perp}\cap M_{D}\subseteq (y_{1},z_{1})^{\perp}\cap M_{D}$ by the definition of $\cC(M_{D})$ in Eqn. \eqref{eqn_codeM}.

Suppose that $(y,z)\in V\setminus\{(0,0)\}$ with $\langle(y,z)^{\perp}\cap M_D\rangle=(y,z)^{\perp}$, we show that $c(y,z)$ is minimal. Take $(y_{1},z_{1})\in V$ with $c(y_{1},z_{1})\preceq c(y,z)$, we have \[(y,z)^{\perp}=\langle(y,z)^{\perp}\cap M_D\rangle\subseteq\langle (y_{1},z_{1})^{\perp}\cap M_{D}\rangle \subseteq (y_{1},z_{1})^{\perp}.\]
If $(y_{1},z_{1})\neq (0,0)$, then $\dim(y,z)^{\perp}=\dim(y_{1},z_{1})^{\perp}$ and so $(y,z)^{\perp}=(y_{1},z_{1})^{\perp}$. If $(y_{1},z_{1})=(0,0)$, then $c(y_{1},z_{1})=0$. Hence $c(y_{1},z_{1})=\mu c(y,z)$ for some $\mu\in\F_{q}$.

Conversely,  we assume that $c(y,z)$ is minimal  for some $(y,z)\in V\setminus\{(0,0)\}$. Suppose to the contrary that $\langle(y,z)^{\perp}\cap M_D\rangle\neq(y,z)^{\perp}$. Then we have $\dim\langle(y,z)^{\perp}\cap M_D\rangle<\dim(y,z)^{\perp}$, i.e., $\dim\langle(y,z)\rangle<\dim\langle(y,z)^{\perp}\cap M_D\rangle^{\perp}.$
Thus, there exists $(y_{1},z_{1})\in\langle(y,z)^{\perp}\cap M_D\rangle^{\perp}$ such that $(y_{1},z_{1})$ is linearly independent with $(y,z)$.  It follows that $(y,z)^{\perp}\cap M_D\subseteq (y_{1},z_{1})^{\perp}\cap M_D$, that is, $c(y_{1},z_{1})\preceq c(y,z)$, which is a contradiction to the minimality of $c(y,z)$. Therefore, we must have $\langle(y,z)^{\perp}\cap M_D\rangle=(y,z)^{\perp}$ as desired if $c(y,z)$ is minimal.

The last claim now follows by the fact that a code is minimal if and only if all codewords of this code are minimal.
\end{IEEEproof}\vspace*{10pt}

\begin{thm}\label{thm_SNC2}
	Take notation as in Theorem \ref{thm_SNC1}. Then $\cC(M_D)$ is a minimal linear code if and only if the following two conditions hold:
	\begin{enumerate}
		\item The set $\overline{D}$ spans $\F_{q^m}$ over $\F_q$, i.e., $\langle \overline{D}\rangle =\F_{q^m}.$
		\item For any $y\in\F_q$ and $z\in\F_{q^m}^*$, $D_{(y,z)}\ne \emptyset$ and $P_{(y,z)}\subseteq \langle z\rangle$, where $D_{(y,z)}$ and $P_{(y,z)}$ are defined in Eqn. \eqref{eqn_Dyz} and Eqn. \eqref{eqn_pyz} respectively.
	\end{enumerate}
\end{thm}
\begin{IEEEproof}
For each codeword $c(y,z)\in\cC(M_D)$ with $(y,z)\in V\setminus\{(0,0)\}$, it is minimal if and only if $\langle(y,z)^{\perp}\cap M_D\rangle=(y,z)^{\perp}$ by Theorem \ref{thm_SNC1}. We split the proof into two cases according as $z=0$ or not.\vspace*{10pt}

\noindent
{\bf Case 1:} In the case $y\in \F_q^*$ and $ z=0$, we compute that
\[(y,0)^{\perp}=\{(u,v)\in V:\,uy+\tr_{\F_{q^m}/\F_{q}}(v\cdot0)=0\}=\{(0,v):v\in\F_{q^m}\}.\]
It follows that
	$(y,0)^{\perp}\cap M_D=\{(0,v): v\in \overline{D}\}.$
Thus, the codeword $c(y,0)$ is minimal if and only if $\langle (y,0)^{\perp}\cap M_D\rangle=(y,0)^{\perp}$, i.e., $\overline{D}$ spans $\F_{q^m}$ over $\F_q$. \vspace*{10pt}
	
\noindent{\bf Case 2:} In the case $(y,z)\in V$ with $z\ne 0$, we have  $(y,z)^{\perp}=\{(u,v) \in V:\, uy+\tr_{\F_{q^m}/\F_{q}}(vz)=0\}.$
 Since $(y,z)^{\perp}$  is $\F_q$-linear and of dimension $m$, it has the following decomposition:
	\begin{equation}\label{eqn_decom}
	 (y,z)^{\perp}=\langle (1,v_0)\rangle \oplus \langle \{(0,v):v\in L(z)\}\rangle,
	\end{equation}
	where $v_0\in\F_{q^m}$ satisfies that $y+\tr_{\F_q^m/\F_q}(v_0z)=0$. In addition, it is clear that
\begin{equation}
\label{eqn_(y,z)^perpcapM}
(y,z)^{\perp}\cap M_D=\{(1,r):r \in D_{(y,z)}\} \cup \{(0,r):r \in \overline{D}_z\}.
\end{equation}

We claim that $\langle (y,z)^{\perp}\cap M_D\rangle=(y,z)^{\perp}$ if and only if $D_{(y,z)}\ne \emptyset$ and $\langle D_{(y,z)} D_{(y,z)}^{(-1)}\cup\overline{D}_z \rangle=L(z)$.
Suppose that $\langle (y,z)^{\perp}\cap M_D\rangle=(y,z)^{\perp}$. Then we have $D_{(y,z)}\ne \emptyset$, otherwise $(y,z)^{\perp}\cap M_D=\{(0,r):r \in \overline{D}_z\}$ and $(1,v_0)\notin\langle (y,z)^{\perp}\cap M_D\rangle$, which contradicts to the decomposition of $(y,z)^{\perp}$ as in Eqn.  \eqref{eqn_decom}. For each $d_0\in D_{(y,z)}$, we have $(1,v_0)-(1,d_0)=(0,v_0-d_0)\in\langle \{(0,v):\,v\in L(z)\}\rangle$ since $\tr_{\F_{q^m}/\F_q}(z(v_0-d_0))=-y+y=0$. Replacing $v_0$ by $d_0$ if necessary, the decomposition  \eqref{eqn_decom} can be replaced as
\begin{equation}
\label{eqn_(y,z)^perp}
(y,z)^{\perp}=\langle (1,d_0)\rangle \oplus \langle \{(0,v):v\in L(z)\}\rangle.
\end{equation}
Since $(1,d_{0})\in(y,z)^{\perp}\cap M_D$, then $\langle (y,z)^{\perp}\cap M_D\rangle=(y,z)^{\perp}$ if and only if
\[\langle\{(0,d_i-d_0):d_i\in D_{(y,z)}\} \cup\{(0,d):d\in \overline{D}_z\} \rangle=\langle\{(0,v):v\in L(z)\} \rangle.\]
This holds for all $d_{0}$ in $D_{(y,z)}$, so we deduce that $\langle D_{(y,z)} D_{(y,z)}^{(-1)}\cup\overline{D}_z \rangle=L(z)$. Conversely, suppose that $D_{(y,z)}\ne \emptyset$ and $\langle D_{(y,z)} D_{(y,z)}^{(-1)}\cup\overline{D}_z \rangle=L(z)$, we take $d_{0}\in D_{(y,z)}$. By above arguments, we replace $v_{0}$ by $d_{0}$ in Eqn. \eqref{eqn_decom} and then get the decomposition \eqref{eqn_(y,z)^perp}. By the condition that $(1,d_{0})\in (y,z)^{\perp}\cap M_{D}$, Eqns. \eqref{eqn_(y,z)^perpcapM}, \eqref{eqn_(y,z)^perp} and $\langle D_{(y,z)} D_{(y,z)}^{(-1)}\cup\overline{D}_z \rangle=L(z)$, we have $\langle (y,z)^{\perp}\cap M_D\rangle=(y,z)^{\perp}$. This proves the claim.

We now show that $\langle D_{(y,z)}D_{(y,z)}^{(-1)}\cup\overline{D}_z \rangle=L(z)$ if and only if $P_{(y,z)}\subseteq \langle z\rangle$. Note that $D_{(y,z)}D_{(y,z)}^{(-1)}\cup\overline{D}_z$ is a subset of $L(z)$. By applying Lemma \ref{lem_1} to $U_1=D_{(y,z)} D_{(y,z)}^{(-1)}\cup\overline{D}_z$ and $V_1=L(z)$, we deduce that $\langle D_{(y,z)}D_{(y,z)}^{(-1)}\cup\overline{D}_z \rangle=L(z)$ if and only if $L(D_{(y,z)}D_{(y,z)}^{(-1)}\cup\overline{D}_z)\subseteq L(L(z))$. The claim follows from the definition of $P_{(y,z)}$ in Eqn. \eqref{eqn_pyz} and $L(L(z))=\langle z\rangle$.

To sum up, we complete this proof by combining above two cases.
\end{IEEEproof}\vspace*{10pt}

Suppose that $D\subseteq\F_{q^m}^{*}$ is $\F_q^*$-invariant. It is straightforward to check that $\overline{D}, \overline{D}_z$ and $ D_{0,z}$ are also $\F_q^*$-invariant for all $z\in \F_{q^m}^*$. By applying Lemma \ref{lem_2} to Eqn. \eqref{eqn_pyz}, we have
\begin{align*}
P_{(y,z)}&=L\left(D_{(y,z)}D_{(y,z)}^{(-1)}\right)\cap L\left(\overline{D}_z\right) \\
&=\{a\in \F_{q^m}:\, \Psi_{a}(\overline{D}_z)=|\overline{D}_z|,\,
\Psi_{a\lambda}\left(D_{(y,z)}D_{(y,z)}^{(-1)}\right)=|D_{(y,z)}|^2,\,\forall\,\lambda\in\mathbb{F}_{q}^{*}\}.
\end{align*}
Since
\begin{align*}
\Psi_{a\lambda}\left(D_{(y,z)}D_{(y,z)}^{(-1)}\right)
=\Psi_{a\lambda}\left(D_{(y,z)}\right)\overline{\Psi_{a\lambda}\left(D_{(y,z)}\right)}
=|\Psi_{a\lambda}\left(D_{(y,z)}\right)|^2,
\end{align*}
we deduce that
\begin{align} 
\label{eqn_pyz2}
P_{(y,z)}=\,\{a\in\F_{q^m}:\,\Psi_{a}(\overline{D}_{z})=|\overline{D}_{z}|,\,
|\Psi_{a\lambda}\left(D_{(y,z)}\right)|=|D_{(y,z)}|,\,\forall\,\lambda\in\mathbb{F}_{q}^{*}\}.
\end{align}

Based on the above arguments and Theorem \ref{thm_SNC2}, we have the following corollary.\vspace*{10pt}

\begin{cor}\label{cor_SNC2}
	Take notation as in Theorem \ref{thm_SNC2} and set $D\subseteq\F_{q^m}^{*}$ to be $\F_q^*$-invariant. Then $\cC(M_D)$ is a minimal linear code if and only if the following two conditions hold:
	\begin{enumerate}
		\item The set $\overline{D}$ spans $\F_{q^m}$ over $\F_q$, i.e., $\langle \overline{D}\rangle =\F_{q^m}.$
		\item For any $y\in\F_q$ and $z\in\F_{q^m}^*$, $D_{(y,z)}\ne \emptyset$ and $P_{(y,z)}\subseteq \langle z\rangle$, where $D_{(y,z)}$ and $P_{(y,z)}$ are defined in Eqn. \eqref{eqn_Dyz} and Eqn. \eqref{eqn_pyz2} respectively.
	\end{enumerate}
\end{cor}
\subsection{Minimal linear codes arising from partial difference sets}
In this subsection, we use $\F_q^*$-invariant partial difference sets to construct minimal linear codes. Suppose that $D\subseteq\F_{q^m}^{*}$ is $\F_q^*$-invariant. Let $\psi$ and $\Psi$ be the canonical additive characters of $\F_q$ and $\F_{q^m}$ respectively. For a property $\cX$, define the Kronecker delta function $[[\cX]] $ as follows:
\begin{equation*}\label{eqn_kro_fun}
   [[\cX]]=\left\lbrace\begin{array}{cl}
  1, & \text{if the property $\cX$ holds,}\\
  0, & \text{otherwise.}
   \end{array}\right.
\end{equation*}
For any $(y,z)\in V$ with $z\ne 0$, we compute the size of the set $D_{(y,z)}$:

\begin{align}
\notag |D_{(y,z)}|&=\sum_{x\in D}[[\tr_{\F_{q^m}/\F_q}(xz)+y=0]]=\frac{1}{q}\sum_{x\in D}\sum_{\lambda\in\mathbb{F}_{q}}\psi(\lambda (\tr_{\F_{q^m}/\F_q}(xz)+y))\\
&=\frac{1}{q}\left(|D|+\sum_{\lambda\in\mathbb{F}_{q}^{*}}\psi(\lambda y)\Psi(\lambda zD)\right)=\left\{\begin{array}{ll}
\frac{1}{q}\left(|D|-\Psi(zD)\right), & \text{ if }y\neq0, \vspace*{1mm} \\
\frac{1}{q}\left(|D|+(q-1)\Psi(zD)\right), & \text{ if }y=0.
\end{array}
\right.\label{eqn_Dyz2}
\end{align}

The last equality holds since $D$ is $\F_q^*$-invariant, and then $\Psi(\lambda zD)=\Psi(zD)$ for any $\lambda\in\F_{q}^{*}$.\vspace*{10pt}

\begin{lemma}\label{lem_pxy}
 Let $D$ be an $\F_q^*$-invariant subset of $\F_{q^m}^*$, and let $P_{(y,z)}$ be defined in Eqn. \eqref{eqn_pyz2} for any fixed $(y,z)\in V$ with $z\in \F_{q^m}^*$. If there exists $a\in P_{(y,z)}$ such that $a\notin \langle z\rangle$, then

 \begin{align}
 	 &|D|=q^{m}+\sum_{\lambda_{1}\in \mathbb{F}_{q}}\Psi\left((\lambda_{1}z+a)D\right)-(q-1)\Psi(zD), \label{eqn_lem_D}\\
 	&|D_{(y,z)}|=\frac{1}{q}|\sum_{\lambda_{1}\in \mathbb{F}_{q}}\psi(\lambda_{1}y)\Psi\left((\lambda_{1}z+a)D\right)|, \label{eqn_lem_Dyz}
 \end{align}
where $\psi$ and $\Psi$ are the canonical additive characters of $\F_q$ and $\F_{q^m}$ respectively.
\end{lemma}
\begin{IEEEproof}
	By Eqn. \eqref{eqn_overlineDz} and the definition of $\overline{D}$,  we have
	\begin{equation}\label{|OverlinDz|}
	|\overline{D}_{z}|=q^{m-1}-1-|D_{(0,z)}|.
	\end{equation}
     We further compute  that
     \begin{align*}
       \Psi(a\overline{D}_{z})=&\sum_{x\in \mathbb{F}_{q^m}^{*}}\Psi(ax)[[\tr_{\F_{q^m}/\F_q}(xz)=0]]-\sum_{x\in D}\Psi(ax)[[\tr_{\F_{q^m}/\F_q}(xz)=0]]\\
        =&\frac{1}{q}\sum_{x\in \mathbb{F}_{q^m}^{*}}\Psi(ax)\sum_{\lambda_{1}\in\mathbb{F}_{q}}\psi\left(\lambda_{1}\tr_{\F_{q^m}/\F_q}(xz)\right) -\frac{1}{q}\sum_{x\in D}\Psi(ax)\sum_{\lambda_{1}\in\mathbb{F}_{q}}\psi\left(\lambda_{1}\tr_{\F_{q^m}/\F_q}(xz)\right)\\
        =&\frac{1}{q}\sum_{\lambda_{1}\in\mathbb{F}_{q}}\sum_{x\in \mathbb{F}_{q^m}^{*}}\Psi\left((a+\lambda_{1} z)x\right)-\frac{1}{q}\sum_{\lambda_{1}\in \mathbb{F}_{q}}\Psi\left((a+\lambda_{1} z)D\right)
     \end{align*}

     Since $a \notin \langle z\rangle$, we have $a+\lambda_{1} z\neq0$ for all $\lambda_{1}\in\mathbb{F}_{q}$. Thus, we have \[\frac{1}{q}\sum_{\lambda_{1}\in\mathbb{F}_{q}}\sum_{x\in \mathbb{F}_{q^m}^{*}}\Psi\left((a+\lambda_{1} z)x\right)=-1\]
      and then
     \begin{equation}
     \label{psi(aODz)}
     \Psi(a\overline{D}_{z})=-1-\frac{1}{q}\sum_{\lambda_{1}\in \mathbb{F}_{q}}\Psi\left((a+\lambda_{1} z)D\right).
     \end{equation}

      By Eqn. \eqref{eqn_pyz2} and  $a\in P_{(y,z)}$, we have $\Psi(a\overline{D}_{z})=|\overline{D}_{z}|$. Then we apply Eqns. \eqref{|OverlinDz|} and  \eqref{psi(aODz)} to get
     \begin{equation}\label{eqn_D0z}
     |D_{(0,z)}|=q^{m-1}+\frac{1}{q}\sum_{\lambda_{1}\in \mathbb{F}_{q}}\Psi\left((a+\lambda_{1} z)D\right).\end{equation}
     Combining Eqn. \eqref{eqn_Dyz2} with Eqn. \eqref{eqn_D0z}, we deduce that Eqn. \eqref{eqn_lem_D} holds.

    Since $a\in P_{(y,z)}$, we have $|D_{(y,z)}|=|\Psi(a\lambda D_{(y,z)})|$ for any $\lambda\in\F_q^{*}$. Without loss of generality, we set $\lambda=1$. It follows that
     \begin{align*}
    |D_{(y,z)}|&=|\Psi(aD_{(y,z)})|=|\sum_{x\in D}\Psi(ax)[[\tr_{\F_{q^m}/\F_{q}}(xz)+y=0]]|\\
    &=|\frac{1}{q}\sum_{\lambda_{1}\in \mathbb{F}_{q}}\sum_{x\in D} \Psi(ax)\psi(\lambda_{1}(\tr_{\F_{q^m}/\F_{q}}(zx)+y))|\\
    &=\frac{1}{q}|\sum_{\lambda_{1}\in \mathbb{F}_{q}}\psi(\lambda_{1}y)\Psi\left((\lambda_{1}z+a)D\right)|.
    \end{align*}

   To sum up, we have completed the proof.
\end{IEEEproof}\vspace*{10pt}

With above preparations, we are now ready to give the construction of the minimal linear codes from $\F_q^*$-invariant partial difference sets.\\

\begin{thm}\label{thm_mainpds}
 Let $D\subseteq \F_{q^m}^*$ be an $\F_q^*$-invariant partial difference set with parameters $(q^m,k,\lambda,\mu)$. Let $\theta_{1},\theta_{2}$ be two restricted eigenvalues of $\text{Cay}(\F_{q^m},D)$ with $\theta_{1}>0>\theta_{2}$, and let $\theta_{0}=\max\{|\theta_{1}|,\,|\theta_{2}|\}$. Suppose that $\cC(M_D)$ is a linear code defined in Eqn. \eqref{eqn_codeM} with the defining set $M_D$ given by Eqn. \eqref{eqn_MD}. Then $\cC(M_D)$ is minimal if the eigenvalues $k$, $\theta_{1}$ and $\theta_{2}$ of $\text{Cay}(\F_{q^m},D)$ satisfy that:
 \begin{enumerate}
 	\item $k-\theta_2\ne q^m $;
 	\item $k>\theta_{1}$ and $k>-(q-1)\theta_{2}$;
 \end{enumerate}	
and one of the following conditions
\begin{enumerate}
    \item[3a)]  $k<q^m+q\theta_{2}-(q-1)\theta_{1}$;
    \item[3b)] $k> \max\{q\theta_{0}+\theta_{1},\, q\theta_{0}-(q-1)\theta_{2}\}$;
    \item[3c)] $q^{m-1}+\theta_{2}-\theta_{1}>\theta_{0}$.
\end{enumerate}
\end{thm}
\begin{IEEEproof}
Recall Section \ref{sec:pds} that $k=\Psi_{0}(D)=|D|$ and $\Psi_{a}(D)\in\{\theta_{1},\,\theta_{2}\}$ for any $a\in\F_{q^m}^{*}$. Suppose that the conditions of the theorem hold. We take three steps to show that the two conditions of Corollary \ref{cor_SNC2} hold, from which we deduce that $\cC(M_D)$ is minimal.

\noindent
{\bf Step 1:} We claim that $\overline{D}$ spans $\F_{q^m}$ over $\F_{q}$. Assume to the contrary that $\langle \overline{D}\rangle \ne  \F_{q^m}$, i.e., $L(\overline{D})\ne \{0\}$. Note that $\overline{D}$ is $\F_q^*$-invariant. For $a\in L(\overline{D})\setminus\{0\}$, we deduce from Lemma \ref{lem_2} that	 \[|\overline{D}|=\Psi_a(\overline{D})=\sum_{x\in\F_{q^m}^*}\Psi(ax)-\Psi(aD)=-1-\Psi(aD).\]
Together with the fact that $|\overline{D}|=q^m-1-|D|=q^m-1-k$, we have $k-\Psi(aD)=q^m$. Since $k-\theta_{1}<k<q^m$, then $\Psi(aD)=\theta_2$, which is a contradiction to the condition 1). Thus, $L(\overline{D})=\{0\}$, i.e., $\overline{D}$ spans $\F_{q^m}$.

\noindent
{\bf Step 2:} We claim that $|D_{(y,z)}|> 0$ for $z\in \F_{q^m}^*$. Since $z\neq0$, then $\Psi(zD)=\theta_{1}$ or $\theta_{2}$. The claim follows from Eqn. \eqref{eqn_Dyz2} and condition 2).
	
\noindent
{\bf Step 3:} We claim that $P_{(y,z)}\subseteq \langle z\rangle$ for any $y\in \F_q,z\in\F_{q^m}^*$, where $P_{(y,z)}$ is defined in  Eqn. \eqref{eqn_pyz2}. We prove by the way of contradiction. Suppose that there exists $a\in P_{(y,z)}\setminus \langle z\rangle$ for some fixed $(y,z)\in V$ with $z\in \F_{q^m}^*$. We just need to show that none of these three conditions 3a), 3b), 3c) holds. Since $a\notin\langle z\rangle$,  we have $\lambda_{1}z+a\ne 0$ for all $\lambda_{1}\in \F_q$, and then $\Psi\left((\lambda_{1}z+a)D\right)$ equals $\theta_{1}$ or $\theta_{2}$. Set $\triangle=\frac{1}{q}\sum_{\lambda_{1}\in \mathbb{F}_{q}}\psi(\lambda_{1}y)\Psi\left((\lambda_{1}z+a)D\right)$. By Eqns. \eqref{eqn_lem_D}, \eqref{eqn_lem_Dyz} in Lemma \ref{lem_pxy} and the triangle inequality, we deduce that
\begin{align}
\label{inequ_D}
&|D|=q^{m}+\sum_{\lambda_{1}\in \mathbb{F}_{q}}\Psi\left((\lambda_{1} z+a)D\right)-(q-1)\Psi(zD) \geq q^m+q\theta_{2}-(q-1)\theta_{1}, \\
\label{inequ_Dyz}
&	 |D_{(y,z)}|=|\triangle|\leq \frac{1}{q}\sum_{\lambda_{1}\in\F_{q}}|\Psi((\lambda_{1}z+a)D)|\leq \theta_{0}.
\end{align}
By the inequality \eqref{inequ_D}, the condition 3a) does not hold.

Next consider Eqn. \eqref{eqn_Dyz2}, we deduce that

\begin{equation}\label{eqn_Deqn}
|D|=\left\{\begin{array}{ll}
q|D_{(y,z)}|+\Psi(zD), & \text{ if }y\neq0, \vspace*{1mm} \\
q|D_{(y,z)}|-(q-1)\Psi(zD), & \text{ if }y=0.
\end{array}
\right.
\end{equation}

We apply the inequality \eqref{inequ_Dyz} to Eqn. \eqref{eqn_Deqn} to deduce that \[k\leq q\theta_{0}+\theta_{1}\text{ or } k\leq q\theta_{0}-(q-1)\theta_{2}\]
according as $y\ne0$ or $y=0$. Thus, the condition 3b) does not hold.

By comparing Eqn. \eqref{eqn_Deqn} and Eqn. \eqref{eqn_lem_D}, we have that
\begin{equation}
\label{eqn_Dyz3}
|D_{(y,z)}|=\left\{\begin{array}{ll}q^{m-1}+\frac{1}{q}\sum_{\lambda\in \mathbb{F}_{q}}\Psi\left((a+\lambda z)D\right)-\Psi(zD), \text{ if }y\neq0, \vspace*{1mm} \\
q^{m-1}+\frac{1}{q}\sum_{\lambda\in \mathbb{F}_{q}}\Psi\left((a+\lambda z)D\right), \text{ if }y=0.\end{array}
\right.
\end{equation}

We apply the inequality \eqref{inequ_Dyz} to Eqn. \eqref{eqn_Dyz3} to deduce that
\[q^{m-1}+\theta_{2}-\theta_{1}\leq \theta_{0}\text{ or }q^{m-1}+\theta_{2}\leq \theta_{0}\]
according as $y\ne0$ or $y=0$. Then the condition 3c) does not hold. This proves the claim.

To sum up, we have established that the two conditions of Corollary \ref{cor_SNC2} hold, and so $\cC(M_{D})$ is a minimal linear code. This completes the proof.
\end{IEEEproof}\vspace*{10pt}

\begin{remark}
\label{remark_georne}
We now remark on 2) of Theorem \ref{thm_mainpds}. To make $\cC(M_D)$ minimal, it is necessary for $D$ to satisfy that $|D_{(y,z)}|\ge 1$ for any $y\in \F_q$ and $z\in \F_{q^m}^*$, cf. Corollary \ref{cor_SNC2}. Hence the condition 2) of Theorem \ref{thm_mainpds} is necessary for $\cC(M_D)$ to be minimal. In particular, since $|D_{(y,z)}|\geq 0$, it suffices to verify that $k\neq\theta_{1}$ and $k\neq-(q-1)\theta_{2}$ for 2) in Theorem \ref{thm_mainpds}.
\end{remark}\vspace*{10pt}

\begin{thm}
\label{thm_weight distri}
Take the same notation as in Theorem \ref{thm_mainpds}. Let $m_1$ and $m_2$ be the corresponding multiplicities of two restricted eigenvalues $\theta_{1}$ and $\theta_{2}$ respectively. Then the weight distribution of $\mathcal{C}(M_{D})$ is listed in Table \ref{Table1}. 
\end{thm}
 \begin{table}[H]
	\centering \caption{Weight distribution of $\mathcal{C}(M_{D})$ with $\F_q^*$-invariant PDS $D$}\label{Table1}
\begin{tabular}{|c|c|}
			\hline
			\textbf{Weight} & \textbf{Frequency} \\
			\hline
			0 & 1 \\
			\hline
			$k$ &$q-1$ \\
			\hline
			$q^m-q^{m-1}$ &$q^{m}-1$ \\
			\hline
			$q^m-q^{m-1}+\theta_{1}$ & $m_1(q-1)$  \\
			\hline
			$q^m-q^{m-1}+\theta_{2}$ & $m_2(q-1)$  \\
			\hline
\end{tabular}
\end{table}	
\begin{IEEEproof}
For each $u\in\mathbb{F}_{q}$, $v\in\mathbb{F}_{q^m}$, we denote the weight of $c(u,v)$ as $\omega_{u,v}$. By Eqns. \eqref{eqn_codeM} and  \eqref{eqn_MD}, we compute that
\begin{align*}	\omega_{u,v}=&(q^m-1)-|\{x\in D:\,u+\tr_{\F_{q^m}/\F_q}(vx)=0\}|-|\{x\in\overline{D}:\,\tr_{\F_{q^m}/\F_q}(vx)=0\}|.
\end{align*}
\noindent{\bf Case 1.} If $u=v=0$, then $\omega_{0,0}=(q^m-1)-(|D|+|\overline{D}|)=0$, and it has frequency 1.

\noindent
{\bf Case 2.} If $u\neq0$ and $v=0$, then \[\omega_{u,0}=(q^m-1)-|\overline{D}|=|D|=k,\]
     and it has frequency  $q-1$.
    	
\noindent
{\bf Case 3.} If $u=0$ and $v\neq0$, then
    	\begin{align*}	 \omega_{0,v}&=(q^m-1)-|\{x\in\mathbb{F}_{q^m}^{*}:\,\tr_{\F_{q^m}/\F_q}(vx)=0\}|\\
    	&=(q^m-1)-(q^{m-1}-1)=q^m-q^{m-1},
    	\end{align*}
    and it has frequency  $q^m-1$.

\noindent
{\bf Case 4.} If $u\neq0$ and $v\neq0$, we compute that
    	\begin{align*}
    	\omega_{u,v}&=(q^m-1)-|D_{(u,v)}|-((q^{m-1}-1)-|D_{(0,v)}|)\\ &=(q^m-q^{m-1})-\frac{1}{q}\left(|D|-\Psi(vD)\right)+\frac{1}{q}\left(|D|+(q-1)\Psi(vD)\right)\\
    	&=q^m-q^{m-1}+\Psi(vD),
    	\end{align*}
   	where $\Psi$ is the canonical additive character of $\F_{q^m}$. Here we used Eqn. \eqref{eqn_Dyz2} in the second equality. Recall that two restricted eigenvalues $\theta_{1}$ and $\theta_{2}$ of $\text{Cay}(\F_{q^m},D)$ have multiplicities $m_1$ and $m_2$ respectively, which have been determined in Lemma \ref{lem_eigenvalue}. By the description of srgs in Section \ref{sec:pds}, we have $\Psi(vD)=\theta_1$ or $\theta_2$ with multiplicity $m_{1}$ or $m_{2}$ respectively, where $v\in \F_{q^m}^*$.  Since $\omega_{u,v}$ is independent with the first coordinate $u$, we deduce that
        $$\omega_{u,v}=\left\{\begin{array}{ll} q^m-q^{m-1}+\theta_{1}, & \text{with\ frequency\ }
        m_{1}(q-1),\\ q^m-q^{m-1}+\theta_{2}, &\text{with\ frequency\ }
        m_{2}(q-1),\end{array}\right. \text{ where } u\in\F_{q}^*,\, v\in\F_{q^m}^*.$$

To sum up, we have determined the weight distribution of $\cC(M_D)$, which is as list in Table \ref{Table1}.
\end{IEEEproof}\vspace*{10pt}

\begin{cor}
\label{cor_weight}
Take the same notation as in Theorem \ref{thm_mainpds}, and assume that $\cC_(M_D)$ is minimal. If  $k\in\{q^m-q^{m-1},\,q^m-q^{m-1}+\theta_{1},\,q^m-q^{m-1}+\theta_{2}\}$, then $\cC(M_{D})$ is a three-weight code; otherwise, $\cC(M_{D})$ is a four-weight code.
\end{cor}
\begin{IEEEproof}
We first show that $q^m-q^{m-1}+\theta_{2}\neq0$. Take $z\in\F_{q^m}^{*}$ such that $\Psi(zD)=\theta_{2}$. By Eqn. \eqref{eqn_Dyz2}, we deduce that $k=q|D_{(y,z)}|+\theta_{2}$ for all $y\in \F_q^*$. By Remark \ref{remark_georne}, the condition 2) of Theorem \ref{thm_mainpds} is necessary for $\cC(M_{D})$ to be minimal, then we have $k>-(q-1)\theta_{2}$. We thus deduce that $-\theta_{2}<|D_{(y,z)}|$ for each $y\in \F_q^*$. By Eqn. \eqref{eqn_Dyz}, we have $|D_{(y,z)}|\leq q^{m-1}$. Since $q^m-q^{m-1}\ge q^{m-1}$, it follows that $q^m-q^{m-1}+\theta_{2}\neq0$ as desired.

By Table \ref{Table1}, we deduce  that $\cC(M_D)$ has exactly three nonzero weights if and only if $k\in\{q^m-q^{m-1},\,q^m-q^{m-1}+\theta_{1},\,q^m-q^{m-1}+\theta_{2}\}$. This completes the proof.
\end{IEEEproof}
\medskip

\begin{cor}
\label{cor_AB}
Suppose that the conditions 1), 2) and at least one of 3a), 3b), 3c) in Theorem \ref{thm_mainpds} hold. If $k\leq (q-1)^2q^{m-2}$ and the characteristic function $f_{D}$ of the set $D$ is not linear, then $\mathcal{C}(M_{D})$ is a $[q^m-1,m+1]$ minimal linear code that does not satisfy the AB condition.
\end{cor}

\begin{IEEEproof}
In the beginning of this section, we have shown that $\mathcal{C}(M_{D})=\mathcal{C}(f_{D})$, where $\mathcal{C}(f_{D})$ is defined in Eqn. (\ref{eqn_CfD}). By Lemma \ref{lem_code_para} and Theorem  \ref{thm_mainpds}, $\mathcal{C}(f_{D})$ is a minimal $[q^m-1,m+1]$ linear code  under our assumptions. Let $\omega_{\text{min}}$ and $\omega_{\text{max}}$ be the minimal and maximal nonzero weights of the code $\cC(M_{D})$ respectively. Since $k\leq (q-1)^2q^{m-2}<q^m-q^{m-1}$, we deduce that  $\omega_{\text{min}}\leq k< q^m-q^{m-1}\leq\omega_{\text{max}}$ by Theorem \ref{thm_weight distri}. It follows that $\frac{\omega_{\text{min}}}{\omega_{\text{max}}}\leq\frac{k}{q^m-q^{m-1}}\leq\frac{q-1}{q}$.
This completes the proof.
\end{IEEEproof}

In the sequel, we focus on partial difference sets of Latin square and negative Latin square type.

\begin{thm}\label{thm_Latin}
Let $q=p^e$ with $p$ prime and $\mathbb{F}_{q^{m}}$ be the finite field with $q^m$ elements such that $m\geq4$, $(m,q)\neq(4,2)$ and 2 divides $em$. Suppose that $D\subsetneq\mathbb{F}_{q^{m}}^{*}$ is an $\F_{q}^{*}$-invariant partial difference set with parameters \[(q^m,r(\sqrt{q^m}-\epsilon),\epsilon \sqrt{q^m}+r^2-3\epsilon r,r^2-\epsilon r),\] which is of Latin square type (resp. negative Latin square type) when $\epsilon=1$ (resp. $\epsilon=-1$). Let $\cC(M_D)$ be the associated linear code as in Eqn. \eqref{eqn_codeM}. The code $\cC(M_D)$ is minimal if one of the following conditions holds:
	\begin{enumerate}
	   \item $\epsilon=1$, $r\neq\sqrt{q^m}$ and $r>1$;
	   \item $\epsilon=-1$, $r\neq\sqrt{q^m}-1$ and $r> \frac{(q-1)\sqrt{q^m}}{\sqrt{q^m}+q}$.
	\end{enumerate}
\end{thm}
\begin{IEEEproof}
By Lemma \ref{lem_eigenvalue}, the eigenvalues of $\text{Cay}(\F_{q^m},D)$ are as follows:
\begin{equation}\label{eqn_eigen}
   \begin{split}
	k&=r(\sqrt{q^m}-\epsilon),\\
 \theta_1&=\frac{1}{2}\left(\epsilon \sqrt{q^m}-2\epsilon r+\sqrt{q^m} \right),\\
  \theta_2&=\frac{1}{2}\left(\epsilon \sqrt{q^m}-2\epsilon r-\sqrt{q^m} \right).
   \end{split}
\end{equation}
Recall that we set $\theta_{0}=\max\{|\theta_{1}|,\,|\theta_{2}|\}$. In this case, we have
\begin{equation}
\label{eqn_theta0}
\theta_{0}=\frac{1}{2}\left(|\sqrt{q^m}-2r|+\sqrt{q^m}\right).
\end{equation}
We now show that the three conditions 1), 2) and 3c) in Theorem \ref{thm_mainpds} are satisfied under our assumptions.\\

\noindent
{\bf Step 1:} We first show that $k-\theta_2\neq q^m$ holds. We compute that
\begin{align*}
k-\theta_2-q^m&=r\sqrt{q^m}+\frac{1}{2}(1-\epsilon)\sqrt{q^m}-q^{m}\\ &=\sqrt{q^m}\left(r+\frac{1}{2}(1-\epsilon)-\sqrt{q^m}\right).
\end{align*}
It is routine to check that it is nonzero under each of the two conditions.\\

\noindent
{\bf Step 2:} We next show that $k>\theta_1$ and $k>-(q-1)\theta_2$. We compute that
\begin{align*}
k-\theta_1&=\left(r-\frac{1}{2}\epsilon-\frac{1}{2}\right)\sqrt{q^m}=\left\lbrace\begin{array}{ll}
        (r-1)\sqrt{q^m}, & \text{ when }\epsilon=1,\\
        r\sqrt{q^m}, & \text{ when }\epsilon=-1.
        \end{array}\right.		
\end{align*}
and
\begin{align*}
k+(q-1)\theta_2&=r(\sqrt{q^m}-\epsilon)+\frac{q-1}{2}(\epsilon\sqrt{q^m}-2\epsilon r-\sqrt{q^m} )\\
&=\left\lbrace \begin{array}{ll}
	     r(\sqrt{q^m}-q), & \text{ when }\epsilon=1,\\
	     r(\sqrt{q^m}+q)-\sqrt{q^m}(q-1), & \text{ when }\epsilon=-1.
\end{array}\right.		
\end{align*}
It is routine to check that both are positive under our assumptions.\\

\noindent
{\bf Step 3:} Finally, we show that $q^{m-1}+\theta_{2}-\theta_{1}-\theta_{0}>0$. By Eqn. \eqref{eqn_eigen}, we deduce that $\theta_{1}-\theta_{2}=\sqrt{q^m}$. Together with Eqn. \eqref{eqn_theta0}, we have
\begin{align*} q^{m-1}+\theta_{2}-\theta_{1}-\theta_{0}&=q^{m-1}-\sqrt{q^m}-\frac{1}{2}(|\sqrt{q^m}-2r|+\sqrt{q^m})\\
&= \left\{\begin{array}{ll}
	q^{m-1}-2\sqrt{q^m}+r, &\text{ if $\sqrt{q^m} \geq 2r$};\\
	q^{m-1}-\sqrt{q^m}-r, & \text{ otherwise. }
	\end{array}\right.
\end{align*}
Since $D$ is a proper subset of $\F_{q^m}^*$, we have $0<k <q^m-1$ and so $0<r< \sqrt{q^m}+\epsilon.$ It follows from $m\geq 4$ that $q^{m-1}-2\sqrt{q^m}+r> 0$.  Since $r< \sqrt{q^m}+\epsilon\leq \sqrt{q^m}+1$, we have 	$q^{m-1}-\sqrt{q^m}-r> q^{m-1}-2\sqrt{q^m}-1>0$ by the assumption that $(m,q)\neq(4,2)$. This establishes 3c) of Theorem \ref{thm_mainpds}.

To sum up, we have proved that the conditions 1), 2) and 3c) in Theorem \ref{thm_mainpds} hold under our assumptions, and so $\cC(M_D)$ is minimal by Theorem \ref{thm_mainpds}.
\end{IEEEproof}

\medskip

\begin{thm}
\label{thm_cycmlc}
Let $q=p^e$ with $p$ prime,  and let  $\mathbb{F}_{q^{m}}$ be a finite field with $q^m$ elements such that $m\geq4$, $(m,q)\neq(4,2)$ and 2 divides $em$. Take the same notation as in Lemma \ref{PDS} and take a proper subset $J\subset\mathbb{Z}_{N}$ of size $u$. If $q$ is odd, we further assume that $N|\frac{q^m-1}{2}$ and $J+\frac{q^m-1}{2}=J(\text{mod }N)$. Let $D=\bigcup_{j\in J}C_{j}$ and $\cC(M_D)$ be the associated linear code defined in Eqn. \eqref{eqn_codeM} with $M_D$ as given by Eqn. \eqref{eqn_MD}. The code $\mathcal{C}(M_{D})$ is a minimal linear code if $J$ is invariant under the map $\rho: \,j\rightarrow j+\frac{q^m-1}{q-1}\pmod N$ and one of the following conditions holds:
\begin{enumerate}
\item $t$ is odd, $u\neq \frac{\sqrt{q^m}N}{\sqrt{q^m}+1}$ and $u>\frac{N}{\sqrt{q^m}+1}$;
\item $t$ is even and $u>\frac{(q-1)\sqrt{q^m}N}{(\sqrt{q^m}+q)(\sqrt{q^m}-1)}$.
\end{enumerate}
\end{thm}
\begin{IEEEproof}
By Lemma \ref{PDS}, the strongly regular Cayley graph $\textup{Cay}(\mathbb{F}_{q^m},D)$ has parameters
\[(q^m,r(\sqrt{q^m}-\epsilon),\epsilon \sqrt{q^m}+r^2-3\epsilon r,r^2-\epsilon r),\,\ \]
where $\epsilon=(-1)^{t+1}$ and $r=\frac{u}{N}\left(\sqrt{q^m}+\epsilon\right)$. It is of Latin square type (resp. negative Latin square type) when $t$ is odd (resp. $t$ is even). Since $J$ is invariant under $\rho$, we deduce that $D$ is $\F_{q}^*$-invariant by Lemma \ref{lem_Fqinv}.  Then the desired result follows from Theorem \ref{thm_Latin}.
\end{IEEEproof}\vspace*{10pt}

\begin{remark}
We consider the special case where $m$ is even and $t$ is odd in Theorem \ref{thm_cycmlc}. In this case, $\sqrt{q^m}=p^{\ell_{1}t}\equiv-1(\text{mod }N)$, so we have $N|(\sqrt{q^m}+1)$. It follows that  $N|\frac{q^m-1}{q-1}$, since $\frac{q^m-1}{q-1}=\frac{\sqrt{q^{m}}-1}{q-1}\cdot(\sqrt{q^{m}}+1)$. We deduce that any subset $J$ of $\mathbb{Z}_{N}$ is trivially invariant under the map $\rho: \,j\rightarrow j+\frac{q^m-1}{q-1}\pmod N$. It follows from Lemma \ref{lem_Fqinv} that $D$ is $\mathbb{F}_{q}^{*}$-invariant. If $q$ is odd, then the conditions $N|\frac{q^m-1}{2}$ and $J+\frac{q^m-1}{2}\equiv J(\text{mod }N)$ in the theorem also hold trivially. Therefore, the requirements on the set $J$ in Theorem \ref{thm_cycmlc} reduce to $|J|\neq \frac{\sqrt{q^m}N}{\sqrt{q^m}+1}$ and $|J|>\frac{N}{\sqrt{q^m}+1}$ in this case.
Since $N|(\sqrt{q^m}+1)$, we have $\frac{N}{\sqrt{q^m}+1}\leq1$, and so the condition $|J|>\frac{N}{\sqrt{q^m}+1}$ holds if $|J|>1$. Therefore, there  are lots of minimal linear codes arising from this construction.
\end{remark}\vspace*{10pt}

 The authors \cite{bonini2019minimal} proposed a useful way to construct minimal linear codes from cutting vectorial $(1,m-1)$-blocking sets in Lemma \ref{lem_cutting}.  But the condition 1) of this lemma is not necessary for a linear code to be minimal. The following example obtained from Theorem \ref{thm_cycmlc} shows that there exists a minimal linear code $\mathcal{C}(f_{D})$ such that $\overline{D}$ is not a cutting vectorial blocking set.\vspace*{10pt}

The next two examples show that the construction in Theorem \ref{thm_cycmlc} is not covered by the cutting vectorial blocking set approach.

\begin{example}
\label{example_cyc}
Let $p=2$, $e=2$, $q=p^e=4$, $m=4$. Set $\gamma$ to be a primitive element of $\mathbb{F}_{4^4}$. Take $\ell_{1}=2$ and $N=p^{\ell_{1}}+1$, then $C_{i}=\{\gamma^{jN+i}:\,0\leq j\leq \frac{q^m-1}{N}-1\}$ for $0\leq i\leq N-1$. We set $J=\mathbb{Z}_{N}\setminus\{0\}$ and $D=D_{J}=\bigcup_{j\in J}C_{j}$. We have $\overline{D}=C_{0}$. Let $\mathcal{C}(M_{D})$ be a linear code defined in Eqn. \eqref{eqn_codeM} with $M_{D}$ defined by Eqn. (\ref{eqn_MD}). We compute that $N=p^{\ell_1}+1=5$ and $t=em/(2\ell_{1})=2$.  Since $N$ divides $\frac{q^m-1}{q-1}$, we deduce from  Lemma \ref{PDS} and Lemma \ref{lem_Fqinv} that $D$ is an $\F_q^*$-invariant  partial difference set. The function $f_{D}$ is not linear by Lemma \ref{lem_fnonlinear}. It is routine to compute that $u=4> \frac{(q-1)\sqrt{q^m}N}{(\sqrt{q^m}+q)(\sqrt{q^m}-1)}$, $k=204$, $\theta_{1}=12$ and $\theta_{2}=-4$.  By Theorem \ref{thm_cycmlc}, $\cC(M_{D})$ is a $[255,5]$ minimal linear code. It is easy to check that $q^{m}-q^{m-1}+\theta_{1}=k$. The code $\cC(M_{D})$ is a three-weight code by Corollary \ref{cor_weight}. Now consider the two affine hyperplanes $H_{1}=\{x\in\mathbb{F}_{4^4}:
\,\tr_{\F_{4^4}/\F_{4}}(x)=0\}$
and $H_{2}=\{x\in\mathbb{F}_{4^4}:\,\tr_{\F_{4^4}/\F_{4}}(\gamma^7x)=0\}$ through the origin.
We compute that $\overline{D} \cap H_{1}=\{1,\gamma^{170},\gamma^{85}\}$ and
 $\overline{D} \cap H_{2}=\{1,\gamma^{170},\gamma^{85},\gamma^{190},\gamma^{20}, \gamma^{105},\gamma^{145},\gamma^{230},\gamma^{60},\gamma^{180},\gamma^{10},\gamma^{95},\gamma^{45},\gamma^{130},\gamma^{215}\}$.
Observe that $\overline{D} \cap H_{1}\subseteq\overline{D} \cap H_{2}$ and so $\overline{D}$ is not a cutting vectorial $(1,3)$-blocking set.
\end{example}\vspace*{10pt}

There are plenty of constructions of $\F_{q}^{*}$-invariant PDS. In Table \ref{Table2}, we list five specific examples taken from Table 1 of \cite{momihara2014lifting}. For more constructions, please refer to \cite{feng2015constructions,feng2012strongly,ge2013constructions, Ma1994, momihara2013strongly, momihara2014lifting}.
\vspace*{10pt}

\begin{example}
\label{example_sporadic}
Let $N$ be a proper divisor of $q^m-1$. Let $\gamma$ be a fixed primitive element of $\F_{q^m}$ and $D=\{\gamma^{iN}:\,0\leq i\leq \frac{q^m-1}{N}-1\}$.  We consider each example in Table \ref{Table2}. It is easy to check that $N|\frac{q^m-1}{q-1}$ in each case, and thus $D$ is $\F_q^*$-invariant by Lemma \ref{lem_Fqinv}. It follows that the characteristic function $f_D$ is not linear by Lemma \ref{lem_fnonlinear}. We can verify that the conditions $1),\,2),\,3a)$ in Theorem \ref{thm_mainpds} hold and so $\cC(M_{D})$ is a $[q^m-1,m+1]$ minimal linear code in each case. We verify that $k\leq (q-1)^2q^{m-2}$ and $k\notin\{q^m-q^{m-1},\,q^m-q^{m-1}+\theta_{1},\,q^m-q^{m-1}+\theta_{2}\}$ in each case, so $\cC(M_D)$ has four nonzero weights and violates the AB condition by Corollary \ref{cor_weight} and Corollary \ref{cor_AB}.

By \cite{momihara2014lifting}, we have $D=\{\gamma^{11i}:\,0\leq i\leq \frac{3^5-1}{11}-1\}$ in the first line of Table \ref{Table2}, where $\gamma$ is a fixed primitive element of $\F_{3^5}$. Set $\overline{D}:=\F_{q^m}^*\setminus D$, which is $\F_{3}^{*}$-invariant. The eigenvalues of $\textup{Cay}(\F_{q^m},\overline{D})$ are $k=220$, $\theta_{1}=4$, $\theta_{2}=-5$. The function $f_{\overline{D}}$ is not linear by Lemma \ref{lem_fnonlinear}. We verify that the conditions $1),\,2),\, 3b)$ of Theorem \ref{thm_mainpds} hold for $\overline{D}$. Hence, $\cC(M_{\overline{D}})$ is a $[3^5-1,6]$ minimal linear code. We check that $k\notin\{q^m-q^{m-1},\,q^m-q^{m-1}+\theta_{1},\,q^m-q^{m-1}+\theta_{2}\}$, and then deduce that $\cC(M_{\overline{D}})$ has four nonzero weights from Corollary \ref{cor_weight}. We have checked by Magma \cite{Magma} that $D$ is not a cutting vectorial $(1,4)$-blocking set, so the code $\cC(M_{\overline{D}})$ does not arise from Lemma \ref{lem_cutting}.

\begin{table}[H]
\centering
\caption{Some examples of PDS $D$'s and the eigenvalues of $\textup{Cay}(\F_{q^m},D)$}\label{Table2}
\begin{small}
	\begin{tabular}{|c|c|c|c|c|c|c|}
			\hline
			No.&$q$ & $m$& $N$& $k$& $\theta_{1}$ & $\theta_{2}$ \\
			\hline
			1&3 & 5&11&22&4&-5 \\
			\hline
			2&5 & 9&19&102796&296&-329 \\
			\hline
			3&3 & 12&35&15184&118&-125 \\
			\hline
			4&7 & 9&37&1090638&584&-1817 \\
			\hline
			5&11 & 7&43&453190&650&-681 \\
			\hline
\end{tabular}
\end{small}
\end{table}	
\end{example}\vspace*{10pt}

\subsection{The automorphism group of $\cC(M_{D})$}
\label{subsec_auto}

Let $D$ be a proper $\F_q^*$-invariant subset of $\F_{q^m}$. An \textit{automorphism} of $D$ is a bijective $\F_q$-linear transformation $g$ of $\F_{q^m}$ that preserves the set $D$, i.e., $\{g(x):\,x\in D\}=D$. We write $\textup{Aut}(D)$ for the set of all automorphisms of $D$.
For each $g\in\text{Aut}(D)$ and $c(u,v)\in\cC(M_{D})$, we define
\[
c(u,v)^{g}:=\left(uf_{D}(g(x))+\tr_{\F_{q^m}/\F_{q}}(vg(x))\right)_{x\in\F_{q^m}^{*}}.
\]
Since $g\in\Aut(D)$, we have $f_D(g(x))=f_D(x)$ for each $x\in\F_{q^m}^{*}$. Recall that $g$ can be regarded as a reduced $q$-polynomial over $\F_{q^m}$, cf.  Section \ref{sec:basicfacts}. Take $\tilde{g}$ to be the trace dual of $g$. It follows that $c(u,v)^{g}=c(u,\tilde{g}(v))$. It is easy to check that $\widetilde{g}$ is also a bijective $\F_q$-linear transformation of $\F_{q^m}$. Thus, $g$ induces an automorphism of the code $\cC(M_{D})$.

In the case $|\Aut(D)|$ is large, the resulting code $\cC(M_{D})$ will also have a large automorphism group. It is thus possible that $\cC(M_{D})$ will have a fast decoding algorithm in such  case. The following is an example where $D$ has a large automorphism group and the associated code $\cC(M_D)$ is minimal.

\begin{example} \label{example_3}
Let $Q:\,\F_{q}^m\rightarrow\F_{q}$ be a nondegenerate quadratic form with $m\ge4$ even, $q=p^h$ and $(m,q)\neq(4,2)$. By Theorem 2.6 of \cite{Ma1994}, $D=\{{\bf x}\in\F_{q}^{m}\setminus\{{\bf 0}\}:\,Q({\bf x})=0\}$ is a PDS with parameters $(q^m,r(\sqrt{q^m}-\epsilon),\epsilon\sqrt{q^m}+r^2-3\epsilon r,r^2-\epsilon r)$. If $Q$ defines a hyperbolic quadric, then $\epsilon=1$, $r=q^{m/2-1}+1$ and $D$ is a PDS of Latin square type. In this case, the automorphism group of $D$ is $\Gamma\textup{O}^{+}(m,q)$ with order $2hq^{m(m-2)/4}(q-1)\Pi_{i=1}^{m/2}(q^{2i}-1)$; if $Q$ defines an elliptic quadric, then $\epsilon=-1$, $r=q^{m/2-1}-1$ and $D$ is a PDS of negative Latin square type. In this case, the  automorphism group of $D$ is $\Gamma\textup{O}^{-}(m,q)$ with order $2hq^{m(m-2)/4}(q-1)(q^{m/2}+1)\Pi_{i=1}^{m/2-1}(q^{2i}-1)$. In both cases, the code $\cC(M_D)$ is minimal by Theorem \ref{thm_Latin}.
\end{example}
\section{Minimal linear codes and secret sharing schemes}
\label{sec:Applications}


In \cite{massey1993minimal} and \cite{massey1995some}, Massey showed that minimal linear codes can be used to construct secret sharing schemes. Later, Yuan and Ding\cite{yuan2005secret} gave a more detailed description on the secret sharing scheme based on a linear code, which we describe below. Let $\cC$ be a $[n,k,d;q]$ linear code with generator matrix  $G=(g_{1},\,g_{2},\,\ldots,\,g_{n})$. The secret $s$ is an element of $\F_{q}$. There are $n-1$ participants $P_{2},\,\ldots,\,P_{n}$ and a trusted person as a dealer in this scheme. The dealer chooses randomly a vector ${\bf u}\in\F_{q}^{k}$ such that $s=t_1={\bf u}g_{1}$, and compute the vector ${\bf t}=(t_{1},\,t_{2},\,\ldots,\,t_{n})={\bf u}G$ and then distribute each $t_{i}$ to participant $P_{i}$ as share for each $i\ge 2$. We call this scheme as the \textit{secret sharing scheme based on $\cC$}. In this scheme, the secret $s=t_1={\bf u}g_{1}$, then a set of shares $\{t_{i_{1}},\,t_{i_{2}},\,\ldots,\,t_{i_{\ell}}\}$ for $2\leq i_{1}<i_{2}<\ldots<i_{\ell}\leq n$ determines this secret if and only if $g_{1}$ is a linear combination of $g_{i_{1}},\,g_{i_{2}},\,\ldots,\,g_{i_{\ell}}$. Thus, a set of participants $P=\{P_{i_{1}},\,P_{i_{2}},\,\ldots,\,P_{i_{\ell}}\}$ for $2\leq i_{1}<i_{2}<\ldots<i_{\ell}\leq n$ determines the secret if and only if there exists a codeword $(1,\,0,\,\ldots,\,0,\,c_{i_{1}},\,0,\,\ldots,\,0,\,c_{i_{\ell}},\,0,\,\ldots,\,0)$ in the dual code $\cC^{\perp}$, cf. \cite[Proposition 1]{yuan2005secret} or \cite{massey1993minimal}. A set of participants is called a \textit{minimal access set} if all participants of the set can recover the secret with their shares but any proper subset can not do so. Clearly, if $\cC^{\perp}$ is a minimal linear code, then there is a one to one correspondence between the set of minimal access sets and the set of the codewords of $\cC^{\perp}$ whose first coordinate is 1. A secret sharing scheme is called democratic of degree $t$ ($t\ge1$) if every group of $t$ participants is in the same number of minimal access sets.

In the sequel, we consider the secret sharing scheme based on $\cC(M_{D})^{\perp}$, where $\cC(M_D)$ is the $[q^m-1,m+1]$ minimal linear code obtained in Section \ref{sec:Mainresults}. Such a secret sharing scheme has many interesting structures, cf. \cite[Proposition 2]{yuan2005secret}. Recall that $\cC(M_{D})=\cC(f_{D})$ and the codeword of $\cC(M_{D})$ is given by $c(u,v)=\left(H(x_{j}\right))_{j=1,\ldots,q^{m}-1}$ for each $(u,v)\in V$, where $H(x)=uf_{D}(x)+\tr_{q^m/q}(vx)$ and $x_1,\cdots,x_{q^m-1}$ is an ordering of the elements of $\F_{q^m}^*$. Let $P_i$ be the participant corresponding to the coordinate labeled by $x_i$, $2\le i\le q^m-1$. By the preceding paragraph, the number of minimal access sets (i.e., the number of minimal codeword in $\cC(M_{D})$ with $H(x_{1})=1$) is $(q^{m+1}-q^m)/(q-1)=q^m$. The number of minimal access sets that the participant $P_{i}$ ($i\ge 2$) lies in is
\[N_{x_{i}}=|\{(u,v)\in V:\,H(x_{1})=1, H(x_{i})\ne0\}|.\]
By basic computation of linear algebra, we get
\[N_{x_{i}}=\left\{\begin{array}{ll}q^m, \text{ if }x_{1}\in \overline{D},\,x_{i}\in\{ax_{1}:\,a\in\F_{q}^{*}\}, \vspace*{1mm} \\
q^m-q^{m-1}, \text{ otherwise. }\end{array}
\right.\]
For instance, in the case that $x_{1}\in \overline{D}$ and $x_{i}\in\{ax_{1}:\,a\in\F_{q}^{*}\}$ for some $2\leq i\leq q^m-1$, we have $N_{x_{i}}=|\{(u,v)\in V:\,\tr_{\F_{q^m}/\F_{q}}(vx_{1})=1\}|=q\cdot q^{m-1}=q^m$. The computation  is similar and technical in other cases, so we omit the details. Actually, this result on $N_{x_{i}}$ is coincide with the Proposition 2 of \cite{yuan2005secret}, and the main step  is to  check in which case $\left(f_{D}(x_1),x_1\right)$ and $\left(f_{D}(x_i),x_i\right)$ are $\F_{q}$-linear dependent.

There are two possible types of schemes arising from $\cC(M_{D})^{\perp}$ according to the choice of $x_1$. If we take $x_{1}\in \overline{D}$, then the corresponding participant $P_{i}$ such that $x_{i}=ax_{1}$ for some $a\in\F_{q}^{*}$ appears in every minimal access set, and such a participant is called \textit{dictatorial}. Such a scheme is useful in scenarios where the bosses must participate in every decision making. If we choose $x_1\in D$, then in this scheme every participant $P_i$, $2\le i\le q^m-1$, lies in $q^m-q^{m-1}$ out of $q^m$ minimal access sets. Such a scheme is democratic of degree at least $1$.

%

\section{Concluding Remarks}
\label{sec:Concluding}
In this paper, we present a general construction of minimal linear codes from $\F_{q}^{*}$-invariant partial difference sets and studied their properties. Our construction yields many minimal linear codes that do not satisfy the AB condition and some examples do not arise from the cutting vectorial blocking set approach. We also show that an automorphism of the partial difference set induces an automorphism of the associated minimal linear code. In the cases where the set has a large automorphism group, the resulting  code will also have a large automorphism group and so potentially have a fast decoding algorithm. Finally, we consider the properties of the secret sharing schemes based on the dual of our codes.

\section*{Acknowledgement}
This work was supported by National Natural Science Foundation of China under Grant No. 11771392. The authors would like to thank the Associate Editor, Prof. Amos Beimel and the reviewers for their constructive comments and suggestions that helped to improve the quality and presentation of this paper.

\bibliographystyle{IEEETrans}
\bibliography{mlcreference}

\end{document}